\documentclass{article}
\usepackage{amsmath,amssymb, amsthm}

\newtheorem{theorem}{Theorem}

\newtheorem{corollary}{Corollary}

\newtheorem{example}{Example}

\newtheorem{lemma}{Lemma}

\newtheorem{remark}{Remark}

\begin{document}

\title{Pointwise summability of Fourier--Laguerre series of integrable
functions}
\author{Bogdan Szal, W\l odzimierz \L enski and Maciej Kubiak \\
University of Zielona G\'{o}ra\\
Faculty of Mathematics, Computer Science and Econometrics\\
65-516 Zielona G\'{o}ra, ul. Szafrana 4a, Poland\\
B.Szal@wmie.uz.zgora.pl,\\
W.Lenski@wmie.uz.zgora.pl,\\
Ma.Kubiak@wmie.uz.zgora.pl}
\date{}
\maketitle

\begin{abstract}
We present an approximation version of the results of D. P. Gupta [ J. of
Approx. Theory, 7 (1973), 226-238] A. N. S. Singroura [Proc. Japan Acad., 39
(4) (1963), 208-210] and G. Szeg\"{o} [Math. Z., 25 (1926), 87-115]. Some
corollaries and examples will also be given.

\ \ \ \ \ \ \ \ \ \ \ \ \ \ \ \ \ \ \ \ 


\ \ \ \ \ \ \ \ \ \ \ \ \ \ \ \ \ \ \ 

\end{abstract}

\section{Introduction}

Let $L_{w}\ $be the class of all real--valued functions, integrable in the
Lebesgue sense over $%
\mathbb{R}
^{+}$ with the weight $w\left( t\right) =e^{-t}t^{\alpha }$ $\left( \alpha
>-1\right) $, i.e.%
\begin{equation*}
\int_{_{_{%
\mathbb{R}
^{+}}}}e^{-t}t^{\alpha }\mid f(t)\mid dt<\infty .
\end{equation*}%
We will consider the Fourier--Laguerre series 
\begin{equation*}
S^{\left( \alpha \right) }f\left( x\right) :=\sum_{\nu =0}^{\infty }a_{\nu
}^{\left( \alpha \right) }(f)L_{\nu }^{\left( \alpha \right) }\left(
x\right) \text{, with }\alpha >-1,
\end{equation*}%
where 
\begin{equation*}
L_{n}^{\left( \alpha \right) }\left( x\right) =\frac{x^{-\alpha }e^{x}}{n!}%
\frac{d^{n}}{dx^{n}}\left( x^{n+\alpha }e^{-x}\right) =\sum_{\nu =0}^{n}%
\frac{\left( -1^{\nu }\right) }{\nu !}\binom{n+\alpha }{n-\nu }x^{\nu }
\end{equation*}%
and%
\begin{equation*}
a_{\nu }^{\left( \alpha \right) }(f)=\frac{1}{\Gamma \left( \alpha +1\right)
A_{n}^{\left( \alpha \right) }}\int_{0}^{\infty }e^{-y}y^{\alpha }L_{\nu
}^{\left( \alpha \right) }\left( y\right) f\left( y\right) dy,\text{ with }%
A_{n}^{\left( \alpha \right) }=\binom{n+\alpha }{n}.
\end{equation*}

Let define the $\left( C,\gamma \right) $ - means$\ $of \ partial sums%
\begin{equation*}
S_{k}^{\left( \alpha \right) }f\left( x\right) =\sum_{\nu =0}^{n}a_{\nu
}^{\left( \alpha \right) }(f)L_{\nu }^{\left( \alpha \right) }\left( x\right)
\end{equation*}%
of $S^{\left( \alpha \right) }f$ as follows 
\begin{equation*}
S_{n}^{\left( \gamma ,\alpha \right) }f\left( x\right) =\frac{1}{%
A_{n}^{\left( \gamma \right) }}\sum_{k=0}^{n}A_{n-k}^{\left( \gamma
-1\right) }S_{k}^{\left( \alpha \right) }f\left( x\right) ,\text{\ \ \ }%
\left( n=0,1,2,...\right)
\end{equation*}%
and let%
\begin{equation*}
\Delta _{x}f\left( t\right) =f\left( x+t\right) -f\left( x\right) .
\end{equation*}

The deviation $S_{n}^{\left( \gamma ,\alpha \right) }f\left( 0\right)
-f\left( 0\right) $\ was examined in the papers \cite{Sin}, \cite{Sz2} and 
\cite{G} as follows:

\textbf{Theorem A.} \cite[Theorem 1]{G} \textit{Let }$f\in L_{w},$\textit{\ }$%
\alpha >-1$\textit{\ and }$\gamma >\alpha +\frac{1}{2}.$\textit{\ If a
function }$f$\textit{\ satisfies the conditions }%
\begin{equation*}
\int_{0}^{u}e^{-t}t^{\alpha }\left\vert \Delta _{0}f\left( t\right)
\right\vert dt=o\left( u^{\alpha +1}\right) \text{ \ \textit{as} }%
u\rightarrow 0^{+}
\end{equation*}%
\textit{and}%
\begin{equation*}
\int_{1}^{\infty }e^{-\frac{t}{2}}t^{\alpha -\gamma -\frac{1}{3}}\left\vert
\Delta _{0}f\left( t\right) \right\vert dt<\infty ,
\end{equation*}%
\textit{then}%
\begin{equation*}
\left\vert S_{n}^{\left( \gamma ,\alpha \right) }f\left( 0\right) -f\left(
0\right) \right\vert =o\left( 1\right) \text{ \textit{as} }n\rightarrow
\infty .
\end{equation*}

Similar results in a case of norm approximation due of C. Markett and E. L.
Poiani in papers \cite{M} and \cite{P} were obtained.

We will say that a nonnegative function $\omega $ is a function of the
modulus of continuity type if it is nondecreasing continuous function on $%
[0,\infty )$ having the following conditions: $\omega \left( 0\right) =0$
and $\omega \left( \delta _{1}+\delta _{2}\right) \leq \omega \left( \delta
_{1}\right) +\omega \left( \delta _{2}\right) $ for any $\delta _{1},\delta
_{2}\in \lbrack 0,\infty )$.

In this paper, we will study the upper bound of the quantity $\left\vert
S_{n}^{\left( \gamma ,\alpha \right) }f\left( 0\right) -f\left( 0\right)
\right\vert $ by some means of a function of the modulus of continuity type $%
\omega .$ From our result we will derive some corollaries, remark and
construct some examples.

\section{Statement of the results}

First we present the estimate of the quantity $\left\vert S_{n}^{\left(
\gamma ,\alpha \right) }f\left( 0\right) -f\left( 0\right) \right\vert $.

\begin{theorem}
Let $f\in L_{w},$ $\alpha >-1$, $\gamma >\alpha +\frac{1}{2}$\ and\ let a
function $\omega $ of the modulus of continuity type satisfy the conditions: 
\begin{equation}
\frac{u^{-\left( \alpha +1\right) }}{\Gamma \left( \alpha +1\right) }%
\int_{0}^{u}e^{-\frac{t}{2}}t^{\alpha }\left\vert \Delta _{0}f\left(
t\right) \right\vert dt=O\left( \omega \left( u\right) \right) \text{ }%
\left( u>0\right)  \label{21}
\end{equation}%
and%
\begin{equation}
\frac{1}{\Gamma \left( \alpha +1\right) }\int_{u}^{\infty }e^{-\frac{t}{2}%
}t^{\alpha -\gamma -\frac{1}{3}}\left\vert \Delta _{0}f\left( t\right)
\right\vert dt=O\left( \omega \left( 1/u\right) \right) \text{ }\left( u\geq
1\right) .  \label{23}
\end{equation}%
Then 
\begin{equation*}
\left\vert S_{n}^{\left( \gamma ,\alpha \right) }f\left( 0\right) -f\left(
0\right) \right\vert =O\left( n^{\eta +\frac{2\left( \alpha -\gamma \right)
+1}{4}}\right) \sum_{k=1}^{n}\frac{\omega \left( 1/k\right) }{k^{^{\eta +%
\frac{2\left( \alpha -\gamma \right) +1}{4}+1}}}+O\left( \omega \left(
1/n^{\eta }\right) \right)
\end{equation*}%
for $0<\eta <-\frac{2\left( \alpha -\gamma \right) +1}{4}.$
\end{theorem}

Now, we formulate some corollaries and remark.

\begin{corollary}
Under the assumptions of the above theorem%
\begin{equation*}
\left\vert S_{n}^{\left( \gamma ,\alpha \right) }f\left( 0\right) -f\left(
0\right) \right\vert =o\left( 1\right) \text{ as }n\rightarrow \infty .
\end{equation*}
\end{corollary}

\begin{remark}
Using Theorem 1 and Corollary 1 we obtain the result of D. P. Gupta from
Theorem A.
\end{remark}

\begin{corollary}
Analyzing the proof of Theorem 1 we can obtain, under the assumptions of
this theorem, the following more precise estimate%
\begin{equation*}
\left\vert S_{n}^{\left( \gamma ,\alpha \right) }f\left( 0\right) -f\left(
0\right) \right\vert =O\left( n^{\frac{2\left( \alpha -\gamma \right) +1}{4}%
}\omega \left( n^{\eta }\right) \right) +O\left( \omega \left( 1/n^{\eta
}\right) \right) +O\left( \omega \left( 1/n\right) \right) ,
\end{equation*}%
when $\frac{2\left( \alpha -\gamma \right) +1}{4}+1<0$.

In the special case, taking $\eta =-\frac{2\left( \alpha -\gamma \right) +1}{%
8}$ , we have%
\begin{equation*}
\left\vert S_{n}^{\left( \gamma ,\alpha \right) }f\left( 0\right) -f\left(
0\right) \right\vert =O\left( \omega \left( 1/n^{\eta }\right) \right) ,
\end{equation*}%
when $\eta \leq 1$.
\end{corollary}

\section{Examples}

Let $f_{1}(t)=e^{-\frac{t}{2}}$ and $\omega _{1}(t)=t$ for $t\geqslant 0$.

It is clear that $f_{1}\in L_{w}$. Moreover, applying the Lagrange mean
value theorem we get that 
\begin{equation*}
|\Delta _{0}f_{1}(t)|=|e^{-\frac{t}{2}}-1|\leqslant \frac{t}{2}
\end{equation*}%
for $t\geq 0$. Therefore, by elementary calculations we get 
\begin{equation*}
\frac{u^{-(\alpha +1)}}{\Gamma (\alpha +1)}\int_{0}^{u}e^{-\frac{t}{2}%
}t^{\alpha }|\Delta _{0}f_{1}(t)|dt
\end{equation*}%
\begin{equation*}
\leqslant \frac{u^{-(\alpha +1)}}{2\Gamma (\alpha +1)}\int_{0}^{u}t^{\alpha
+1}dt=\frac{1}{2\Gamma (\alpha +1)(\alpha +2)}\omega _{1}(u)
\end{equation*}%
for $u>0$ and 
\begin{equation*}
\frac{1}{\omega _{1}\left( \frac{1}{u}\right) \Gamma (\alpha +1)}%
\int_{u}^{\infty }e^{-\frac{t}{2}}t^{\alpha -\gamma -\frac{1}{3}}|\Delta
_{0}f_{1}(t)|dt
\end{equation*}%
\begin{equation*}
\leqslant \frac{u}{2\Gamma (\alpha +1)}\int_{u}^{\infty }e^{-\frac{t}{2}%
}t^{\alpha -\gamma +\frac{2}{3}}dt\leqslant \frac{1}{2\Gamma (\alpha +1)}%
\int_{u}^{\infty }e^{-\frac{t}{2}}t^{\alpha -\gamma +\frac{5}{3}}dt
\end{equation*}%
\begin{equation*}
\leqslant \frac{1}{2\Gamma (\alpha +1)}\int_{0}^{\infty }e^{-\frac{t}{2}%
}t^{\alpha -\gamma +\frac{5}{3}}dt=\frac{1}{2\Gamma (\alpha +1)}%
\int_{0}^{\infty }e^{-\frac{t}{2}}t^{-1+\alpha -\gamma +\frac{8}{3}}dt
\end{equation*}%
\begin{equation*}
=\frac{1}{\Gamma (\alpha +1)}2^{\alpha -\gamma +\frac{5}{3}}\Gamma \left(
\alpha -\gamma +\frac{8}{3}\right) <\infty
\end{equation*}%
for $u\geq 1$ and $\alpha -\gamma +\frac{8}{3}>0$.

Hence the function $f_{1}$ satisfies the conditions $(1)$ and $(2)$. Using
Theorem 1 we get the following estimate for $|S_{n}^{(\gamma ,\alpha
)}f_{1}(0)-f_{1}(0)|$:

\begin{example}
Let $\alpha >-1,$ $\alpha +\frac{1}{2}<\gamma <$ $\alpha +\frac{8}{3}$ and $%
0<\eta <-\frac{2(\alpha -\gamma )+1}{4}$. Then 
\begin{equation*}
|S_{n}^{(\gamma ,\alpha )}f_{1}(0)-f_{1}(0)|=O\left( n^{\eta +\frac{2(\alpha
-\gamma )+1}{4}}\right) \sum\limits_{k=1}^{n}\frac{1}{k^{\eta +\frac{%
2(\alpha -\gamma )+1}{4}+2}}+O\left( \frac{1}{n^{\eta }}\right) .
\end{equation*}
\end{example}

Suppose $f_{2}\left( t\right) =t^{\delta }$ and $\omega _{2}\left( t\right)
=t^{\delta }$ for $\delta \in (0,1]$ and $t\geqslant 0$.

Obviously $f_{2}\in L_{w}$. In addition, it is easy to show that%
\begin{equation*}
\frac{u^{-(\alpha +1)}}{\Gamma (\alpha +1)}\int_{0}^{u}e^{-\frac{t}{2}%
}t^{\alpha }|\Delta _{0}f_{2}(t)|dt
\end{equation*}%
\begin{equation*}
\leqslant \frac{u^{-(\alpha +1)}}{\Gamma (\alpha +1)}\int_{0}^{u}t^{\alpha
+\delta }dt=\frac{1}{\Gamma (\alpha +1)(\alpha +\delta +1)}\omega _{2}(u)
\end{equation*}%
for $u>0$ and \ 
\begin{equation*}
\frac{1}{\omega _{2}\left( \frac{1}{u}\right) \Gamma (\alpha +1)}%
\int_{u}^{\infty }e^{-\frac{t}{2}}t^{\alpha -\gamma -\frac{1}{3}}|\Delta
_{0}f_{2}(t)|dt
\end{equation*}%
\begin{equation*}
\leqslant \frac{u^{\delta }}{\Gamma (\alpha +1)}\int_{u}^{\infty }e^{-\frac{t%
}{2}}t^{\alpha -\gamma +\delta -\frac{1}{3}}dt\leqslant \frac{1}{\Gamma
(\alpha +1)}\int_{u}^{\infty }e^{-\frac{t}{2}}t^{\alpha -\gamma +2\delta -%
\frac{1}{3}}dt
\end{equation*}%
\begin{equation*}
\leqslant \frac{1}{\Gamma (\alpha +1)}\int_{0}^{\infty }e^{-\frac{t}{2}%
}t^{\alpha -\gamma +2\delta -\frac{1}{3}}dt=\frac{1}{\Gamma (\alpha +1)}%
\int_{0}^{\infty }e^{-\frac{t}{2}}t^{-1+\alpha -\gamma +2\delta +\frac{2}{3}%
}dt
\end{equation*}%
\begin{equation*}
=\frac{1}{\Gamma (\alpha +1)}2^{\alpha -\gamma +2\delta +\frac{2}{3}}\Gamma
\left( \alpha -\gamma +2\delta +\frac{2}{3}\right) <\infty
\end{equation*}%
for $u\geq 1$ and $\alpha -\gamma +2\delta +\frac{2}{3}>0$.

Therefore the function $f_{2}$ satisfies the conditions $(1)$ and $(2)$.
Using Theorem 1 we get the following estimate for $|S_{n}^{(\gamma ,\alpha
)}f_{2}(0)-f_{2}(0)|$:

\begin{example}
Let $\alpha >-1,$ $\delta \in (0,1]$, $\alpha +\frac{1}{2}<\gamma <\alpha
+2\delta +\frac{2}{3}$ and $0<\eta <-\frac{2(\alpha -\gamma )+1}{4}$. Then 
\begin{equation*}
|S_{n}^{(\gamma ,\alpha )}f_{2}(0)-f_{2}(0)|=O\left( n^{\eta +\frac{2(\alpha
-\gamma )+1}{4}}\right) \sum\limits_{k=1}^{n}\frac{1}{k^{\eta +\frac{%
2(\alpha -\gamma )+1}{4}+1+\delta }}+O\left( \frac{1}{n^{\eta \delta }}%
\right) .
\end{equation*}
\end{example}

\section{Auxiliary results}

We begin\ this section by some notations from \cite{Sz} . We have%
\begin{equation*}
L_{k}^{\left( \alpha +1\right) }\left( y\right) =\sum_{\nu =0}^{k}L_{\nu
}^{\left( \alpha \right) }\left( y\right) ,\text{ \ \ }L_{\nu }^{\left(
\alpha \right) }\left( 0\right) =\binom{\nu +\alpha }{\nu }
\end{equation*}%
and therefore 
\begin{eqnarray*}
S_{k}^{\left( \alpha \right) }f\left( 0\right) &=&\frac{1}{\Gamma \left(
\alpha +1\right) }\int_{0}^{\infty }e^{-y}y^{\alpha }L_{k}^{\left( \alpha
+1\right) }\left( y\right) f\left( y\right) dy, \\
S_{n}^{\left( \gamma ,\alpha \right) }f\left( 0\right) &=&\frac{1}{\Gamma
\left( \alpha +1\right) A_{n}^{\left( \gamma \right) }}\int_{0}^{\infty
}e^{-y}y^{\alpha }L_{n}^{\left( \alpha +\gamma +1\right) }\left( y\right)
f\left( y\right) dy.
\end{eqnarray*}%
Hence, by evidence equality%
\begin{equation*}
\frac{1}{\Gamma \left( \alpha +1\right) }\int_{0}^{\infty }e^{-y}y^{\alpha
}L_{\nu }^{\left( \alpha +1\right) }\left( y\right) dy=\left\{ 
\begin{array}{c}
1\text{\ \ if }\nu =0, \\ 
0\text{ \ if }\nu \neq 0,%
\end{array}%
\right.
\end{equation*}%
we have%
\begin{eqnarray*}
S_{k}^{\left( \alpha \right) }f\left( 0\right) -f\left( 0\right) &=&\frac{1}{%
\Gamma \left( \alpha +1\right) }\int_{0}^{\infty }e^{-y}y^{\alpha
}L_{k}^{\left( \alpha +1\right) }\left( y\right) \Delta _{0}f\left( y\right)
dy, \\
.S_{n}^{\left( \gamma ,\alpha \right) }f\left( 0\right) -f\left( 0\right) &=&%
\frac{1}{\Gamma \left( \alpha +1\right) A_{n}^{\left( \gamma \right) }}%
\int_{0}^{\infty }e^{-y}y^{\alpha }L_{n}^{\left( \alpha +\gamma +1\right)
}\left( y\right) \Delta _{0}f\left( y\right) dy
\end{eqnarray*}

Next, we present the useful estimates:

\begin{lemma}
\cite[p. 172]{Sz} Let $\beta $ be an arbitrary real number, $c$ and $\delta $
be fixed positive constants. Then%
\begin{equation*}
\left\vert L_{n}^{\left( \beta \right) }\left( x\right) \right\vert =\left\{ 
\begin{array}{c}
O\left( n^{\beta }\right) \text{ \ \ \ \ \ \ \ \ \ \ \ \ \ \ \ \ \ \ \ \ \ \
\ \ \ if }0\leq x\leq \frac{c}{n}, \\ 
O\left( x^{-\left( 2\beta +1\right) /4}n^{\left( 2\beta -1\right) /4}\right) 
\text{ \ \ if }\frac{c}{n}\leq x\leq \delta .%
\end{array}%
\right.
\end{equation*}
\end{lemma}

\begin{lemma}
\cite[p. 235]{Sz} Let $\beta $ and $\lambda $ be arbitrary real numbers, $%
\delta >0$ and $0<\theta <4$ . Then%
\begin{equation*}
\underset{x}{\max }e^{-x/2}x^{\lambda }\left\vert L_{n}^{\left( \beta
\right) }\left( x\right) \right\vert =\left\{ 
\begin{array}{c}
O\left( n^{\max \left( \lambda -\frac{1}{2},\frac{\beta }{2}-\frac{1}{4}%
\right) }\right) \text{ \ \ \ if }\delta \leq x\leq \left( 4-\theta \right)
n, \\ 
O\left( n^{\max \left( \lambda -\frac{1}{3},\frac{\beta }{2}-\frac{1}{4}%
\right) }\right) \text{ \ \ if \ \ \ \ \ \ \ \ \ \ \ \ \ \ \ \ \ \ }x\geq
\delta .%
\end{array}%
\right.
\end{equation*}
\end{lemma}

\begin{lemma}
\cite[Vol. I, (1.15) and Theorem 1.17]{AZ} If $\gamma >-1,$ then 
\begin{equation*}
A_{n}^{\left( \gamma \right) }=\binom{n+\gamma }{n}\simeq O\left( \left(
n+1\right) ^{\gamma }\right) \text{ }
\end{equation*}%
and $A_{n}^{\left( \gamma \right) }$ is positive for $\gamma >-1$ increasing
(as a function of $n$) for $\gamma >0$ and decreasing for $-1<\gamma <0.$
\end{lemma}

\section{Proofs of Theorems}

\subsection{Proofs of Theorem 1}

It is clear that if 
\begin{eqnarray*}
&&S_{n}^{\left( \gamma ,\alpha \right) }f\left( 0\right) -f\left( 0\right) \\
&=&\frac{1}{\Gamma \left( \alpha +1\right) A_{n}^{\left( \gamma \right) }}%
\int_{0}^{\infty }e^{-y}y^{\alpha }L_{n}^{\left( \alpha +\gamma +1\right)
}\left( y\right) \Delta _{0}f\left( y\right) dy \\
&=&\left( \int_{0}^{1/n}+\int_{1/n}^{1}+\int_{1}^{n^{\eta }}+\int_{n^{\eta
}}^{\infty }\right) =J_{1}+J_{2}+J_{3}+J_{4}\text{,}
\end{eqnarray*}%
then%
\begin{equation*}
\left\vert S_{n}^{\left( \gamma ,\alpha \right) }f\left( 0\right) -f\left(
0\right) \right\vert \leq \left\vert J_{1}\right\vert +\left\vert
J_{2}\right\vert +\left\vert J_{3}\right\vert +\left\vert J_{4}\right\vert .
\end{equation*}%
By Lemma 1, Lemma 3 and (\ref{21})%
\begin{eqnarray*}
\left\vert J_{1}\right\vert &\leq &\frac{1}{\Gamma \left( \alpha +1\right)
A_{n}^{\left( \gamma \right) }}\int_{0}^{1/n}e^{-\frac{y}{2}}y^{\alpha
}\left\vert L_{n}^{\left( \alpha +\gamma +1\right) }\left( y\right)
\right\vert \left\vert \Delta _{0}f\left( y\right) \right\vert dy \\
&=&\frac{O\left( n^{\alpha +\gamma +1}\right) }{\Gamma \left( \alpha
+1\right) A_{n}^{\left( \gamma \right) }}\int_{0}^{1/n}e^{-\frac{y}{2}%
}y^{\alpha }\left\vert \Delta _{0}f\left( y\right) \right\vert dy \\
&=&\frac{O\left( n^{\alpha +1}\right) }{\Gamma \left( \alpha +1\right) }%
\int_{0}^{1/n}e^{-\frac{y}{2}}y^{\alpha }\left\vert \Delta _{0}f\left(
y\right) \right\vert dy=O\left( \omega \left( 1/n\right) \right) \\
&\leq &O\left( n^{\left( 2\alpha -2\gamma +1\right) /4}\right) \sum_{k=1}^{n}%
\frac{\omega \left( 1/k\right) }{k^{\left( 2\alpha -2\gamma +1\right) /4+1}}
\\
&\leq &O\left( n^{\eta +\frac{2\left( \alpha -\gamma \right) +1}{4}}\right)
\sum_{k=1}^{n}\frac{\omega \left( 1/k\right) }{k^{^{\eta +\frac{2\left(
\alpha -\gamma \right) +1}{4}+1}}}\text{,}
\end{eqnarray*}%
\ with $0<\eta <-\frac{2\left( \alpha -\gamma \right) +1}{4}$.

Using Lemma 1, we get\ 
\begin{eqnarray*}
\left\vert J_{2}\right\vert &\leq &\frac{1}{\Gamma \left( \alpha +1\right)
A_{n}^{\left( \gamma \right) }}\int_{1/n}^{1}e^{-y}y^{\alpha }\left\vert
L_{n}^{\left( \alpha +\gamma +1\right) }\left( y\right) \right\vert
\left\vert \Delta _{0}f\left( y\right) \right\vert dy \\
&=&\frac{O\left( n^{\frac{2\left( \alpha +\gamma \right) +1}{4}}\right) }{%
\Gamma \left( \alpha +1\right) A_{n}^{\left( \gamma \right) }}%
\int_{1/n}^{1}e^{-\frac{y}{2}}y^{\alpha }\left\vert \Delta _{0}f\left(
y\right) \right\vert y^{-\frac{2\left( \alpha +\gamma \right) +3}{4}}dy.
\end{eqnarray*}%
Let $F_{\alpha }\left( y\right) =\frac{y^{-\left( \alpha +1\right) }}{\Gamma
\left( \alpha +1\right) }\int_{0}^{y}e^{-\frac{u}{2}}u^{\alpha }\left\vert
\Delta _{0}f\left( u\right) \right\vert du$. Applying Lemma 3 and
integrating by parts with $\gamma >\alpha +\frac{1}{2}$ and $\alpha >-1$ we
have%
\begin{eqnarray*}
\left\vert J_{2}\right\vert &=&\frac{O\left( n^{\left( 2\alpha +2\gamma
+1\right) /4}\right) }{\Gamma \left( \alpha +1\right) A_{n}^{\left( \gamma
\right) }}\int_{1/n}^{1}e^{-y/2}y^{\left( 2\alpha -2\gamma -3\right)
/4}\left\vert \Delta _{0}f\left( y\right) \right\vert dy \\
&=&O\left( n^{\left( 2\alpha -2\gamma +1\right) /4}\right) \left\{ \left[
F_{\alpha }\left( y\right) y^{\frac{2\left( \alpha -\gamma \right) +1}{4}}%
\right] _{y=1/n}^{1}\right. \\
&&+\left. \frac{2\left( \alpha +\gamma \right) +3}{4}\int_{1/n}^{1}F_{\alpha
}\left( y\right) y^{\frac{2\left( \alpha -\gamma \right) +1}{4}-1}dy\right\}
\\
&\leq &O\left( n^{\left( 2\alpha -2\gamma +1\right) /4}\right) \left\{
F_{\alpha }\left( 1\right) +\int_{1/n}^{1}F_{\alpha }\left( y\right) y^{%
\frac{2\left( \alpha -\gamma \right) +1}{4}-1}dy\right\}
\end{eqnarray*}%
\begin{eqnarray*}
&=&O\left( n^{\left( 2\alpha -2\gamma +1\right) /4}\right) \left\{ F_{\alpha
}\left( 1\right) +\int_{1}^{n}F_{\alpha }\left( 1/y\right) y^{-\frac{2\left(
\alpha -\gamma \right) +1}{4}-1}dy\right\} \\
&=&O\left( n^{\left( 2\alpha -2\gamma +1\right) /4}\right) \left\{ F_{\alpha
}\left( 1\right) +\sum_{k=1}^{n-1}\int_{k}^{k+1}F_{\alpha }\left( 1/y\right)
y^{-\frac{2\left( \alpha -\gamma \right) +1}{4}-1}dy\right\}
\end{eqnarray*}%
\begin{eqnarray*}
&\leq &O\left( n^{\left( 2\alpha -2\gamma +1\right) /4}\right) \left\{
F_{\alpha }\left( 1\right) +\sum_{k=1}^{n-1}F_{\alpha }\left( 1/k\right) k^{-%
\frac{2\left( \alpha -\gamma \right) +1}{4}-1}\right\} \\
&\leq &O\left( n^{\left( 2\alpha -2\gamma +1\right) /4}\right)
\sum_{k=1}^{n-1}2F_{\alpha }\left( 1/k\right) k^{-\frac{2\left( \alpha
-\gamma \right) +1}{4}-1} \\
&\leq &O\left( n^{\left( 2\alpha -2\gamma +1\right) /4}\right) \sum_{k=1}^{n}%
\frac{F_{\alpha }\left( 1/k\right) }{k^{\left( 2\alpha -2\gamma +1\right)
/4+1}}.
\end{eqnarray*}%
By (\ref{21}) we obtain%
\begin{eqnarray*}
\left\vert J_{2}\right\vert &=&O\left( n^{\left( 2\alpha -2\gamma +1\right)
/4}\right) \sum_{k=1}^{n}\frac{\omega \left( 1/k\right) }{k^{\left( 2\alpha
-2\gamma +1\right) /4+1}} \\
&\leq &O\left( n^{\eta +\frac{2\left( \alpha -\gamma \right) +1}{4}}\right)
\sum_{k=1}^{n}\frac{\omega \left( 1/k\right) }{k^{^{\eta +\frac{2\left(
\alpha -\gamma \right) +1}{4}+1}}},
\end{eqnarray*}%
\ with $0<\eta <-\frac{2\left( \alpha -\gamma \right) +1}{4}.$

Applying Lemma 2 \ with $\alpha +\gamma +1$ instead of $\beta $ and $\lambda
=\frac{2\alpha +2\gamma +3}{4}$ (since $\max \left( \lambda -\frac{1}{2},%
\frac{\alpha +\gamma +1}{2}-\frac{1}{4}\right) =\frac{2\alpha +2\gamma +1}{4}
$) we have%
\begin{eqnarray*}
&&\left\vert J_{3}\right\vert \\
&\leq &\frac{1}{\Gamma \left( \alpha +1\right) A_{n}^{\left( \gamma \right) }%
}\int_{1}^{n^{\eta }}e^{-y/2}y^{\left( 2\alpha -2\gamma -3\right)
/4}\left\vert \Delta _{0}f\left( y\right) \right\vert e^{-y/2}y^{\left(
2\alpha +2\gamma +3\right) /4}\left\vert L_{n}^{\left( \alpha +\gamma
+1\right) }\left( y\right) \right\vert dy \\
&=&\frac{O\left( n^{\left( 2\alpha +2\gamma +1\right) /4}\right) }{\Gamma
\left( \alpha +1\right) A_{n}^{\left( \gamma \right) }}\int_{1}^{n^{\eta
}}e^{-y/2}y^{\left( 2\alpha -2\gamma -3\right) /4}\left\vert \Delta
_{0}f\left( y\right) \right\vert dy.
\end{eqnarray*}%
Using Lemma 3 and integrating by parts with $\gamma >\alpha +\frac{1}{2}$
and $\alpha >-1$, we get%
\begin{eqnarray*}
\left\vert J_{3}\right\vert &=&\frac{O\left( n^{\left( 2\alpha -2\gamma
+1\right) /4}\right) }{\Gamma \left( \alpha +1\right) }\int_{1}^{n^{\eta }}%
\left[ e^{-y/2}y^{\alpha }\left\vert \Delta _{0}f\left( y\right) \right\vert %
\right] y^{\left( -2\alpha -2\gamma -3\right) /4}dy \\
&=&O\left( n^{\left( 2\alpha -2\gamma +1\right) /4}\right) \int_{1}^{n^{\eta
}}\frac{d}{dy}\left[ \int_{0}^{y}\frac{e^{-u/2}u^{\alpha }\left\vert \Delta
_{0}f\left( u\right) \right\vert }{\Gamma \left( \alpha +1\right) }du\right]
y^{\left( -2\alpha -2\gamma -3\right) /4}dy
\end{eqnarray*}%
\begin{eqnarray*}
&=&O\left( n^{\left( 2\alpha -2\gamma +1\right) /4}\right) \left\{ \left[
y^{\left( -2\alpha -2\gamma -3\right) /4}\int_{0}^{y}\frac{e^{-u/2}u^{\alpha
}\left\vert \Delta _{0}f\left( u\right) \right\vert }{\Gamma \left( \alpha
+1\right) }du\right] _{1}^{n^{\eta }}\right. \\
&&+\left. \frac{2\alpha +2\gamma +3}{4}\int_{1}^{n^{\eta }}\frac{d}{dy}\left[
\int_{0}^{y}\frac{e^{-u/2}u^{\alpha }\left\vert \Delta _{0}f\left( u\right)
\right\vert }{\Gamma \left( \alpha +1\right) }du\right] y^{\frac{-2\alpha
-2\gamma -3}{4}-1}dy\right\}
\end{eqnarray*}%
\begin{equation*}
\leq O\left( n^{\left( 2\alpha -2\gamma +1\right) /4}\right) \left\{
F_{\alpha }\left( n^{\eta }\right) n^{\eta \frac{2\left( \alpha -\gamma
\right) +1}{4}}+\frac{2\alpha +2\gamma +3}{4}\int_{1}^{n^{\eta }}F_{\alpha
}\left( y\right) y^{\frac{2\left( \alpha -\gamma \right) +1}{4}-1}dy\right\}
.\text{ }
\end{equation*}%
By (\ref{21}) we obtain%
\begin{equation*}
\left\vert J_{3}\right\vert \leq O\left( n^{\left( 2\alpha -2\gamma
+1\right) /4}\right) \left\{ \omega \left( n^{\eta }\right) n^{\eta \frac{%
2\left( \alpha -\gamma \right) +1}{4}}+\frac{2\alpha +2\gamma +3}{4}%
\int_{1}^{n^{\eta }}\omega \left( y\right) y^{\frac{2\left( \alpha -\gamma
\right) +1}{4}-1}dy\right\}
\end{equation*}%
\begin{equation*}
\leq O\left( 1\right) \left\{ n^{\left( 2\alpha -2\gamma +1\right)
/4}n^{\eta }\omega \left( 1\right) n^{\eta \frac{2\left( \alpha -\gamma
\right) +1}{4}}+n^{\left( 2\alpha -2\gamma +1\right) /4}\omega \left(
n^{\eta }\right) \int_{1}^{n^{\eta }}y^{\frac{2\left( \alpha -\gamma \right)
+1}{4}-1}dy\right\}
\end{equation*}%
\begin{eqnarray*}
&\leq &O\left( n^{\left( 2\alpha -2\gamma +1\right) /4}n^{\eta }n^{\eta 
\frac{2\left( \alpha -\gamma \right) +1}{4}}+n^{\frac{2\left( \alpha -\gamma
\right) +1}{4}}n^{\eta }\right) \omega \left( 1\right) \text{ \ }\leq
O\left( n^{\frac{2\left( \alpha -\gamma \right) +1}{4}}n^{\eta }\right)
\omega \left( 1\right) \\
&\leq &O\left( n^{\eta +\frac{2\left( \alpha -\gamma \right) +1}{4}}\right)
\sum_{k=1}^{n}\frac{\omega \left( 1/k\right) }{k^{^{\eta +\frac{2\left(
\alpha -\gamma \right) +1}{4}+1}}},
\end{eqnarray*}%
with $0<\eta <-\frac{2\left( \alpha -\gamma \right) +1}{4}.$

If $\lambda =\gamma +\frac{1}{3}$ then $\lambda -\frac{1}{3}=\gamma >\frac{%
\alpha +\gamma +1}{2}-\frac{1}{4}$ since $\gamma >\alpha +\frac{1}{2}$. So,
applying Lemma 2 with $\alpha +\gamma +1$ instead of $\beta $ and $\lambda
=\gamma +\frac{1}{3}$ we obtain%
\begin{eqnarray*}
&&\left\vert J_{4}\right\vert \\
&\leq &\frac{1}{\Gamma \left( \alpha +1\right) A_{n}^{\left( \gamma \right) }%
}\int_{n^{\eta }}^{\infty }e^{-y/2}y^{\left( 3\alpha -3\gamma -1\right)
/3}\left\vert \Delta _{0}f\left( y\right) \right\vert e^{-y/2}y^{\left(
3\gamma +1\right) /3}\left\vert L_{n}^{\left( \alpha +\gamma +1\right)
}\left( y\right) \right\vert dy \\
&=&\frac{O\left( n^{\gamma }\right) }{\Gamma \left( \alpha +1\right)
A_{n}^{\left( \gamma \right) }}\int_{n^{\eta }}^{\infty }e^{-y/2}y^{\left(
3\alpha -3\gamma -1\right) /3}\left\vert \Delta _{0}f\left( y\right)
\right\vert dy.
\end{eqnarray*}%
Next, using Lemma 3 and (\ref{23}) we get%
\begin{equation*}
\left\vert J_{4}\right\vert =\frac{O\left( 1\right) }{\Gamma \left( \alpha
+1\right) }\int_{n^{\eta }}^{\infty }e^{-y/2}y^{\left( 3\alpha -3\gamma
-1\right) /3}\left\vert \Delta _{0}f\left( y\right) \right\vert dy=O\left(
\omega \left( 1/n^{\eta }\right) \right).
\end{equation*}

Finally, collecting the above estimates we have 
\begin{equation*}
\left\vert S_{n}^{\left( \gamma ,\alpha \right) }f\left( 0\right) -f\left(
0\right) \right\vert =O\left( n^{\eta +\frac{2\left( \alpha -\gamma \right)
+1}{4}}\right) \sum_{k=1}^{n}\frac{\omega \left( 1/k\right) }{k^{^{\eta +%
\frac{2\left( \alpha -\gamma \right) +1}{4}+1}}}+O\left( \omega \left(
1/n^{\eta }\right) \right) .
\end{equation*}%
and\ our proof is completed. $\blacksquare $

\section{Conclusions}

We investigated pointwise approximation of real--valued functions,
integrable in the Lebesgue sense over $%
\mathbb{R}
^{+}$ with the weight $w\left( t\right) =e^{-t}t^{\alpha }$ $\left( \alpha
>-1\right) $\ by the $\left( C,\gamma \right) $ - means$\ $of \ partial sums
of their Fourier--Laguerre series. In particular, we estimated the deviation 
$\left\vert S_{n}^{\left( \gamma ,\alpha \right) }f\left( 0\right) -f\left(
0\right) \right\vert $\ by means of a function of the modulus of continuity
type $\omega $. From our result some corollaries were derived and some
examples were constructed.

\end{document}